\documentclass[11pt,bezier,epsf]{article}
\usepackage{amsmath,amssymb,amsfonts}
\usepackage{graphicx}

 \numberwithin{equation}{section}

\title{On Adaptive Multiple-Shooting Method for Stochastic Multi-Point Boundary Value Problems}
\author{Ali Foroush Bastani\footnote{Corresponding author}\hspace{1cm} Davood Damircheli}
\date{Department of Mathematics, Inistitute for Advanced Studies in Basic Sciences,\\
P.O. Box 45195-1159,  Zanjan, Iran}

\begin{document}
\maketitle

\begin{quotation}
\noindent {\bf Abstract} This paper presents an adaptive
multiple-shooting method to solve stochastic multi-point boundary
value problems. The heuristic to choose the shooting points is
based on separating the effects of drift and diffusion terms and comparing the corresponding solution components with a
pre-specified initial approximation. Having obtained the mesh
points, we solve the underlying stochastic differential
equation on each shooting interval with a first-order strongly-convergent stochastic Runge-Kutta method. We illustrate the effectiveness of this
approach on 1-dimentional and 2-dimentional test problems and compare our results with other non-adaptive alternative techniques proposed in the literature.
\end{quotation}



 \maketitle

\section{Introduction}

Numerical methods for solving initial value problems in
stochastic differential equations (SDE-IVPs) have been extensively
researched in the last two decades (see e.g. \cite{KP,M} and the references therein). This is not the
case for stochastic boundary value problems (SDE-BVPs or SBVPs
for short), because of complications both in theoretical as well
as computational aspects. These equations appear naturally in a
variety of fields such as smoothing \cite{P}, maximum a
posteriori estimation of trajectories of diffusions
\cite{ZD1}, wave motion in random media \cite{FM}, stochastic
optimal control \cite{ZD}, valuation of boundary-linked assets
\cite{EV} and in the study of reciprocal processes \cite{K}. They also arise from the semi-discretization in space of stochastic partial differential equations by the method of lines approximation \cite{MD}.
Taking into account the fact that the exact solution of these
equations are rarely available in analytic form, trying to find
efficient approximation schemes for the trajectories of the
solution process or its moments, seems to be a natural candidate.
During the last years, several authors have studied with different
techniques, the numerical solution of SBVPs of the form:
\begin{equation}\label{SBVP}
\left\{
\begin{array}{lcl}
dX(t)=f(X(t),t)dt+g(X(t),t)\circ dW(t),\quad X(t)\in{\mathbb{R}^d}, \quad 0\leq t\leq T, \\
\alpha(X)=c,\\
\end{array}
\right.
\end{equation}
in which $f:\mathbb{R}^d\times [0,T]\rightarrow \mathbb{R}^d$ and
$g:\mathbb{R}^d \times [0,T]\rightarrow {\mathbb{R}}^{d \times
d}$ are continuous globally Lipschitz functions with polynomial
growth, $W(t)$ is a $d$-dimensional Wiener process,
$\alpha:C^0(\mathbb{R}^d\times [0,T])\rightarrow \mathbb{R}^d$ is
a continuous operator and $c\in\mathbb{R}^d$ is a constant vector.
The existence and uniqueness of the solution process as well as the {\bf Markov
field property} of it have been studied by some authors, among them
we mention \cite{OP,OP1,NP,ZD,G}. Due to the anticipative nature of the
solution process, the main machinery in the study of these
equations have turned out to be the Malliavin calculus \cite{Nu}.

The majority of research in this field has concentrated around
two-point SBVPs (TP-SBVPs) corresponding to the choice
\begin{equation}
\alpha(X)=h(X(0),X(T))=c,
\end{equation}
in which $h:\mathbb{R}^d\times \mathbb{R}^d\rightarrow
\mathbb{R}^d$ is a given (possibly nonlinear) function and $c$ is defined as before. In this category, we must point out to {\bf linear TP-SBVP}s in which both the drift and diffusion coefficients ($f$ and $g$ respectively in (\ref{SBVP})) are linear functions of their arguments and the function $h$ is of the form
\begin{equation}
h(y,z)=H_0y+H_1z,
\end{equation}
in which $H_0$ and $H_1$ are $d\times d$ matrices. At the same time, the special
class of {\bf functional boundary conditions} of the form
\begin{equation}\label{FBC}
\alpha(X)=\int_{0}^{T}dA(t)X(t)=c,
\end{equation}
have also been of interest, in which $A(t)$ is a $d\times d$ matrix
valued integrator. The other interesting case is the multi-point
SBVP (or MP-SBVP for short) having the boundary condition
\begin{equation}\label{MPBC}
\alpha(X)=\sum_{j=1}^{N_s}A_jX(\tau_j)=c,
\end{equation}
in which $A_1,A_2,\cdots,A_m$ are constant square matrices of
order $d$ and $\tau_1,\tau_2,\cdots,\tau_{N_s}\in[0,T]$ are given {\it
switching points} with the property $\tau_i<\tau_j$, for $i<j$.
This boundary condition could be considered as the result of a
quadrature formula applied to approximate the general form
(\ref{FBC}) and will be of special interest in this paper.

On the numerical side, some efforts have been directed towards
devising efficient numerical schemes for (\ref{SBVP}) among them we mention the following: Allen and
Nunn \cite{AN} propose two methods for linear two-dimensional
second order SBVPs, one based on finite differences and the other
based on simple-shooting. They analyze the convergence
properties of these methods and report some numerical
experiments confirming their theoretical results. Arciniega and
Allen \cite{AA} examine a shooting-type method for systems of
linear SBVPs of the form (\ref{SBVP}). This method could be viewed
as a generalization of the {\it complementary function approach} for
deterministic BVPs adopted to solve SBVPs \cite{RS}. Arciniega
\cite{Arc} extends this work to the nonlinear case and performs
some error analysis for this new scheme. Ferrante, Kohatsu-Higa and Sanz-Sol$\acute{e}$
\cite{FKS} use a strong Euler-Maruyama approximation to find
strong solutions of (\ref{SBVP}) with linear boundary conditions.
They obtain error estimates for this method without accompanying any numerical
results to their theoretical findings. In a recent paper, Esteban-Bravo and Vidal-Sanz \cite{EV} use the wavelet-collocation scheme to find approximations to trajectories of the solution for a general version of (\ref{SBVP}) with boundary conditions of the form (\ref{FBC}).
We must also mention the work of Prigarin and Winkler \cite{PW} in which they propose a special member of the general Markov chain Monte Carlo (MCMC) approach namely the Gibbs sampler to construct realizations of the solution process. The convergence is proved for the special case of linear TP-SBVPs and some guidelines have been provided to cope with the general nonlinear case and also boundary value problems for stochastic partial differential equations.

Among the above-mentioned schemes, the simple-shooting method
which relies on transforming the SBVP (\ref{SBVP}) to an SDE-IVP, has
shown to have good accuracy properties, but it may give
unacceptable approximate solutions on long time intervals. This
is specially the case when the underlying SDE is
\textit{unstable} i.e. almost all sample paths are rapidly
growing in absolute value. Our aim here is to circumvent this
deficiency by developing an adaptive multiple-shooting method to
solve (\ref{SBVP}) based on a detailed analysis of the sample
paths of the corresponding stochastic equation. The idea is to adaptively
subdivide the typical interval $[\tau_i,\tau_{i+1}]$ into a grid of {\it
shooting points}
$$\tau_i=t_{i,1}(\omega)<t_{i,2}(\omega)<\cdots<t_{i,j}(\omega)<\cdots<t_{i,N(i)}(\omega)=\tau_{i+1}$$
in which the $t_{i,j}$'s and also $N(i)$ will depend on the
particular realization (indexed by $\omega$) of the underlying Wiener process. In each
interval $[\tau_i,\tau_{i+1}]$, starting from $t_{i,1} = \tau_i$,
the criterion we choose to obtain $t_{i,j+1}$ from $t_{i,j}$ is to use an idea adopted from the {\bf operator-splitting} method
to investigate the behavior of the two local SDE-IVPs arising from the drift and diffusion components of the underlying SDE and controlling
upon their growth on this subinterval. For this purpose, we employ an initial approximation to the solution which (approximately) satisfies the boundary conditions and compare it with the two corresponding SDE-IVP solutions. To obtain the mesh
points, we solve the above mentioned SDE-IVPs on each shooting
interval with a first-order strongly convergent {\bf stochastic Runge-Kutta} method introduced in \cite{BH}. We show that this strategy significantly
enhances the accuracy and stability properties of the simple-shooting method and at the
same time reduces the computational cost of the long-time integration problem to a great extent. Comparison with other schemes like simple-shooting, finite-differences, wavelet-collocation and also the fixed-step multiple shooting method itself, confirms that the proposed method is a reliable alternative than the widely used non-adaptive approaches in the literature.

The rest of this paper is organized as follows. In section 2, we
present the multiple shooting framework to solve SBVPs with
multi-point boundary conditions. The criterion to select the
shooting points which forms the foundation of our adaptive
strategy will be discussed in section 3. The details of optimal
parameter tuning for the proposed scheme and implementation details will be described in section 4. We conclude the paper by commenting on some possible ways to extend this work into more general frameworks.

\section{Multiple Shooting Method for MP-SBVPs}
In this section, we describe the multiple-shooting framework to
approximate the sample paths of the equation (\ref{SBVP}). This
can be considered as the extension of methods presented in
\cite{AA,Arc} and will serve as the ground base for our adaptive
scheme. For this purpose, consider the following MP-SBVP in
Stratonovich form:
\begin{equation}\label{MPSBVP}
\left\{
\begin{array}{lcl}
dX(t)=f(X(t),t)dt+g(X(t),t)\circ dW(t),\quad X(t)\in{\mathbb{R}^d}, \quad 0\leq t\leq T, \\
\sum_{j=1}^{N_s}A_jX(\tau_j)=c.\\
\end{array}
\right.
\end{equation}
Without loss of generality, we assume throughout the paper that
$\tau_1=0$ and $\tau_{N_s}=T$. For each realization of the Wiener
process, we are interested in finding the corresponding realization of the
solution process satisfying (\ref{MPSBVP}). Assume that for each $i=1,2,\cdots,N_s-1$,
$I_i=[\tau_i,\tau_{i+1}]$ is subdivided into the shooting
intervals $[t_{i,j},t_{i,j+1}],~j=1,\cdots,N(i)-1$ with
$t_{i,1}=\tau_{i}$ and $t_{i,N(i)}=\tau_{i+1}$. The adaptive
procedure used to obtain them will be discussed in section 3 but in the sequel, we assume that they are known. If $X_{i,j}(t;s_{i,j})$ solves the local SDE-IVPs:
\begin{equation}\label{Local}
\left\{
\begin{array}{lcl}
dX_{i,j}(t;s_{i,j})=f(X_{i,j}(t;s_{i,j}),t)dt+g(X_{i,j}(t;s_{i,j}),t)\circ
dW(t),\quad t\in[t_{i,j},t_{i,j+1}],\\
X_{i,j}(t_{i,j};s_{i,j})=s_{i,j},
\end{array}
\right.
\end{equation}
for $i=1,2,\cdots,N_s-1$ and $j=1,\cdots,N(i)-1$, augmentation of the local solutions $X_{i,j}(t;s_{i,j})$ with imposition of continuity condition at the interior shooting points and satisfaction of multi-point boundary
conditions at switching points will result in
a global approximation to $X(t)$. For this purpose, we
find the unknown initial conditions $s_{i,j}$'s by solving
the system of $D=d\times (\sum_{i=1}^{N_s-1}[N(i)-1])+1$ nonlinear
equations:
\begin{equation}\label{F}
F({\bf s})=0,
\end{equation}
in which
$${\bf s}=(s_{1,1}^T,\cdots,s_{1,N(1)-1}^T,s_{2,1}^T,\cdots,s_{2,N(2)-1}^T,\cdots,s_{N_s-1,1}^T,\cdots,s_{N_s-1,N(N_s-1)-1}^T,s_{N_s,1}^T)^T\in\mathbb{R}^{D},$$
is the {\it shooting vector} and $F(\textbf{s})$ is given by:
\begin{equation}\label{Nonlinear}
\begin{array}{l}
 F({\bf s}) = \left[ {\begin{array}{cc}
   {{\bf s}}_2  - {\bf X}_1({\bf s}_1) \\
    \vdots   \\
   {{\bf s}}_{N_s}  - {\bf X}_{N_s-1}({\bf s}_{N_s-1})  \\
   {g({\bf s},c)}  \\
\end{array}}\right], \\
   \end{array}
\end{equation}
in which ${\bf s}_j = (s_{j,2},~s_{j,3},\cdots,s_{j,N(j)-1},~s_{j+1,1})^T$ for $j=2,3,\cdots,N_s$,
\begin{equation}
{\bf X}_j = \left[ {\begin{array}{*{20}c}
   X_{j,1}(t_{j,2};s_{j,1})  \\
   \\
   X_{j,2}(t_{j,3};s_{j,2})\\
   \\
    \vdots   \\
   X_{j,N(j)-2}(t_{j,N(j)-1};s_{j,N(j)-2})  \\
   \\
   X_{j,N(j)-1}(t_{j+1,1};s_{j,N(j)-1})  \\
\end{array}} \right],
\end{equation}
for $j=1,\cdots,N_s-1$ and
\begin{equation}
g({\bf s},c) =
A_1s_{1,1}+A_2s_{2,1}+\cdots+A_{N_s-1}s_{N_s-1,1}+A_{N_s}s_{N_s,1}-c.
\end{equation}
The solution of system (\ref{Nonlinear}), which provides a global
refinement of the solution values at the gridpoints, is usually done within the framework of a damped-Newton iteration whose $k$-th iteration
${\bf s}^k$ is of the form
\begin{equation}
{\bf s}^{k+1}={\bf s}^k-\lambda_k[DF({\bf s}^k)]^{-1}F({\bf s}^k).
\end{equation}
In this relation, $\lambda_k\in(0,1]$ is the {\it relaxation} or {\it damping factor} and
$DF({\bf s}^k)$ is the Jacobian matrix of $F({\bf s})$ evaluated
at the $k$-th iteration.

It can be shown that
$$DF({\bf s})=\left[
    \begin{array}{ccccc}
      - {\bf \Gamma}_1 & \textbf{I}_1 &  &   &  \\
       & -{\bf \Gamma}_2 & \textbf{I}_2 &  &  \\
       &  & \ddots &\ddots  &  \\
       &  &  & -{\bf \Gamma}_{N_s-1} & \textbf{I}_{N_s-1} \\
            A_1 & A_2  &   \ldots  & A_{N_s-1} & A_{N_s} \\
    \end{array}
  \right]$$
  in which
\begin{equation}
{\bf \Gamma}_j  = \left[ {\begin{array}{*{20}c}
    \Gamma_{j,1}& {} & {} & \\
    & \Gamma_{j,2} & {}&  \\
   {} & &  \ddots&  \\
   {} & {} & &\Gamma_{j,N(j)}  \\
\end{array}} \right]
\end{equation}
and the components $\Gamma_{j,k}\equiv D_{s_{j,k}}X_{j,k-1}(t_{j,k}
;s_{j,k-1})$ for each $j$ and $k$ are $d\times d$ matrices and
\begin{equation}
\textbf{I}_j = \left[ {\begin{array}{*{20}c}
    I_{d\times d}& {} & {} & \\
    & I_{d\times d} & {}&  \\
   {} & &  \ddots&  \\
   {} & {} & &I_{d\times d} \\
\end{array}} \right]
\end{equation}
is an $N(j)\times N(j)$ identity matrix. It is obvious that the
exact computation of $\Gamma_{j,k}$ requires the analytic solution of the
local SDE-IVPs (\ref{Local}). It is worth pointing out here that
although it is possible to approximate $\Gamma_{j,k}$'s by linearization of
the corresponding local SDEs and integrating them up to
$t_{i,j}$, we will adopt an alternative strategy by approximating
the derivative terms by finite differences. A strategy for
choosing the $\lambda_k$'s has also been developed and thoroughly
tested in \cite{D} which will be pursued here.
\section{Adaptive Sequential Selection of Shooting Points}
Multiple-shooting method as a natural generalization of the simple-shooting idea, significantly enhances the stability properties
of its ancestor and behaves much better than it in terms of accuracy
and rate of convergence and so has been a preferred
choice to solve deterministic boundary value problems
\cite{AMR,Keller,SB}. The main drawback of this method could be
attributed to its computational cost which is directly
proportional to the number of shooting points in the integration
interval. To reduce these costs, some authors have proposed to
devise a control mechanism on the number and location of shooting
points in such a way that the stability and accuracy properties
of the method are preserved. This strategy has the additional
advantage of resolving the special features of the solution in
the integration interval: ``$\dots$ a multiple-shooting approach should
permit step sizes to be chosen sequentially, fine in the boundary
layers, and coarse in the smooth regions'' \cite{EM}. We extend
this argument to the case of non-smooth solutions - the feature
which is intrinsic for SDEs - and show that the adaptive
selection of shooting points based on the driving force for this
non-smooth behavior, i.e. the underlying Wiener process and also comparing the solution with an initial approximate solution, will
have an overall performance much better than the corresponding
fixed step-size counterpart.

To start the adaptive procedure, we first find a simple piecewise
linear approximation to the solution, $\theta(t)$, which
approximately satisfies the multi-point boundary conditions (\ref{MPBC}). To
find this approximation, we discretize the SDE-IVP in each
interval with the Euler-Maruyama method and then solve the following
system of nonlinear equations for $\overline{\theta}_{\tau_j}$'s:
\begin{equation}\label{NLEQ}
\left\{
\begin{array}{lcl}
\overline{\theta}_{\tau_2}=\overline{\theta}_{\tau_1}+(\tau_2-\tau_1)f(\overline{\theta}_{\tau_1},\tau_1)+(W(\tau_2)-W(\tau_1))g(\overline{\theta}_{\tau_1},\tau_1),\\
\overline{\theta}_{\tau_3}=\overline{\theta}_{\tau_2}+(\tau_3-\tau_2)f(\overline{\theta}_{\tau_2},\tau_2)+(W(\tau_3)-W(\tau_2))g(\overline{\theta}_{\tau_2},\tau_2),\\
~\vdots \quad  \quad   \vdots  \\
\overline{\theta}_{\tau_{N_s}}=\overline{\theta}_{\tau_{N_s-1}}+(\tau_{N_s}-\tau_{N_s-1})f(\overline{\theta}_{\tau_{N_s-1}},\tau_{N_s-1})+(W(\tau_{N_s})-W(\tau_{N_s-1}))g(\overline{\theta}_{\tau_{N_s-1}},\tau_{N_s-1}),\\
\\
\sum_{j=1}^{N_s} A_j\overline{\theta}_{\tau_{j}}=c.
\end{array}
\right.
\end{equation}
The continuous piecewise linear approximation could then be obtained by linear interpolation:
\begin{equation}\label{interp1}
\theta(t) =
\frac{t-\tau_i}{\tau_{i+1}-\tau_i}\overline{\theta}_{\tau_{i+1}}+
\frac{\tau_{i+1}-t}{\tau_{i+1}-\tau_i}\overline{\theta}_{\tau_{i}},\quad
t\in [\tau_i,\tau_{i+1}],~~ i =1,\cdots,N_s.
\end{equation}
Consider now the interval $[\tau_i,\tau_{i+1}]$ and put
$t_{i,1}:=\tau_i$. Starting from $t_{i,j}$ and to obtain
the next shooting point in this interval, we integrate the following two local
SDE-IVP's:
\begin{equation}\label{SDE-IVP-1}
\left\{
\begin{array}{lcl}
d\widehat{X}(t)&=&f(\widehat{X}(t),t)dt,\quad \quad t \in [t_{i,j},\tau_{i+1}],\\
\widehat{X}(t_{i,j})&=&\theta(t_{i,j}),\\
\end{array}
\right.
\end{equation}
\begin{equation}\label{SDE-IVP-2}
\left\{
\begin{array}{lcl}
d\widetilde{X}(t)&=&g(\widetilde{X}(t),t)\circ dW(t),\quad \quad t \in [t_{i,j},\tau_{i+1}],\\
\widetilde{X}(t_{i,j})&=&\theta(t_{i,j}),\\
\end{array}
\right.
\end{equation}
by deterministic and stochastic components of an SRK method, described in the next section. We will terminate the integration when we
reach the first point in our discretization satisfying:
\begin{equation}\label{Stop-Loss}
t_{i,j}\leq s \leq \tau_{i+1},\quad ||\widehat{X}(s)||\geq
L_1(s)\quad\textmd{or}\quad||\widetilde{X}(s)||\geq L_2(s),
\end{equation}
in which $L_1(s)$ and $L_2(s)$ will be specified in the sequel. We then put $t_{i,j+1}:=s$ and restart the integration of both (\ref{SDE-IVP-1}) and (\ref{SDE-IVP-2}) from $t_{i,j+1}$ using $\widehat{X}(t_{i,j+1})=\widetilde{X}(t_{i,j+1})=\theta(t_{i,j+1})$ as the initial guess. This procedure will be continued up until the point $t_{i,N(i)}=\tau_{i+1}$ is reached and then will be continued from the next shooting interval to finally arrive at $T$.

The first untold story in our description of the algorithm is the selection of the ``stop-loss functions'' $L_1(s)$ and $L_2(s)$ which control upon the location of our shooting points. The most intuitionistic proposal could be
\begin{equation}\label{Criteria}
L_1(s)=\alpha\|\theta(s)\|,\quad L_2(s)=\beta\|\theta(s)\|
\end{equation}
for some positive constants $\alpha$ and $\beta$ (see e.g. \cite{SB} Section 7.3.6 for a similar idea in the case of deterministic BVPs). We could choose the $\alpha$ and $\beta$ coefficients in (\ref{Criteria}) time-dependent and find an empirical optimal relation for them, but our numerical experience shows that the gain in efficiency is not substantial. We have also tested other stopping criteria based only on the size of the increments of the Wiener process which has resulted in the selection of more shooting points but has not improved the accuracy in a comprehendible way. Another proposal is to find the first point $s$ which simultaneously maximizes the following quantities:
\begin{eqnarray*}\label{New-Criteria}
{\Bbb P}\Big(\frac{\|\widehat{X}(s)\|}{\|\theta(s)\|}\geq\alpha\Big),\quad {\Bbb P}\Big(\frac{\|\widetilde{X}(s)\|}{\|\theta(s)\|}\geq\beta\Big)
\end{eqnarray*}
for given positive $\alpha$ and $\beta$. The idea has led us to solve simple constrained stochastic programming problems in each step (with exact solutions for linear SBVPs) that needs further investigation and will be pursued in a forthcoming paper.

It is evident from the form of our adaptation criteria that in the case of weak driving noise process, we are controlling upon the size of the solution process and look at the first time at which the norm of the solution starts to deviate from the initial piecewise linear approximation. On the other hand, when the increments of the Brownian noise become large in some portions of the solution domain, we must finish the integration and select the current point as a suitable shooting point. In both of these scenarios, we must come back to the initial approximation $\theta(t)$ and continue the integration from the initial value $\theta(t_{i,j+1})$.

Having described the way in which we choose our shooting points for each realization, we now need to tell the other story about the time marching procedure to solve our SDE-IVPs resulting from the multiple shooting method, which will be discussed in the next section.
\subsection{Stochastic Runge-Kutta Family}
Among the many possible choices of methods to integrate the ODE-IVP and SDE-IVP problems in
(\ref{SDE-IVP-1}) and (\ref{SDE-IVP-2}) (assuming w.l.o.g. that both equations are autonomous), we choose to work with a special member from the general class of stochastic Runge-Kutta (SRK) methods of the form
\begin{equation} \left\{
\begin{array}{lcl}
\overline{\eta}_i &=& \overline{X}_{n}+ h\sum_{j=1}^{s}a_{ij}f(\overline{\eta}_j)+J_1\sum_{j=1}^{s}b_{ij}g(\overline{\eta}_j),\quad i=1,\dots,s\\
\\
\overline{X}_{{n+1}}&=&\overline{X}_{n}+h\sum_{j=1}^{s}\alpha_{j}f(\overline{\eta}_j)+J_1\sum_{j=1}^{s}\gamma_{j}g(\overline{\eta}_j),\\
\end{array}
\right.
\end{equation}
in which $\overline{X}_{{n}}$ and $\overline{X}_{{n+1}}$ are
approximations to $X(t_n)$ and $X(t_{n+1})$ respectively,
$h=t_{n+1}-t_n$ and $J_1=\int_{t_n}^{t_{n+1}}\circ
dW(s)=W(t_{n+1})-W(t_n)$. Here, $A=(a_{ij})_{i,j=1}^n$ and $B=(b_{ij})_{i,j=1}^n$ are
$s\times s$ matrices with real elements and
$\alpha^T=(\alpha_1,\dots,\alpha_s)$ and
$\gamma^T=(\gamma_1,\dots,\gamma_s)$ are row vectors in
$\mathbb{R}^s$. A typical member of this family could be
represented by the Butcher tableau
\begin{center}
\begin{tabular}{c|c|c}
& A & B\\
\hline \vspace{-0.3cm}
&  & \\
& $\alpha^T$ & $\gamma^T$ \\
\end{tabular}
\end{center}
and according to the theory presented in \cite{BB}, the highest possible order of strong (and also weak) convergence among all consistent choices for $A,B,\alpha$ and $\gamma$ is one (see \cite{KP} for notions of strong and weak convergence in the SDE literature).
In this work, we use a three stage SRK method (dubbed R3) as
the underlying numerical integrator which has the tableau
\begin{center}
\begin{tabular}{c|c c c|c c c}
& 0 & 0 & 0 & 0 & 0 & 0 \\
& $\tfrac12$ & 0 & 0 & $\tfrac{1}{2}$ & 0 & 0 \\
& 0 & $\tfrac34$ & 0 & 0 & $\tfrac{3}{4}$ & 0 \\
\hline
& $\tfrac29$ & $\tfrac39$ & $\tfrac49$ & $\tfrac29$ & $\tfrac39$ & $\tfrac49$\\
\end{tabular}
\end{center}
and its deterministic and stochastic components are themselves valid numerical integration schemes \cite{BH}. More specifically, we integrate the IVP (\ref{SDE-IVP-1}) with a method of the form
\begin{equation} \left\{
\begin{array}{lcl}
\overline{\eta}_1 &=& \overline{X}_{n},\\
\\
\overline{\eta}_2 &=& \overline{X}_{n} + \frac{h}{2} f(\overline{\eta}_1),\\
\\
\overline{\eta}_3 &=& \overline{X}_{n} + \frac{3h}{4} f(\overline{\eta}_2),\\
\\
\overline{X}_{{n+1}}&=&\overline{X}_{n}+\frac{h}{9}(2f(\overline{\eta}_1)+3f(\overline{\eta}_2)+4f(\overline{\eta}_3)),
\end{array}
\right.
\end{equation}
and integrate the SDE (\ref{SDE-IVP-2}) with another method having the form
\begin{equation} \left\{
\begin{array}{lcl}
\overline{\eta}_1 &=& \overline{X}_{n},\\
\\
\overline{\eta}_2 &=& \overline{X}_{n} + \frac{J_1}{2} g(\overline{\eta}_1),\\
\\
\overline{\eta}_3 &=& \overline{X}_{n} + \frac{3J_1}{4} g(\overline{\eta}_2),\\
\\
\overline{X}_{{n+1}}&=&\overline{X}_{n}+\frac{J_1}{9}(2g(\overline{\eta}_1)+3g(\overline{\eta}_2)+4g(\overline{\eta}_3)).
\end{array}
\right.
\end{equation}
It is interesting to note here that the first scheme has third-order of convergence for a deterministic IVP and this will result in higher precision when we are faced with an SDE-BVP having a weak driving noise. On the other hand and for the second scheme, we have first-order of strong convergence for drift-free SDEs and when the drift is going to diminish in some portions of the problem domain, we have an exact-enough method to trace the non-smooth path of the corresponding realization.
\section{Numerical Experiments}
In this section, we report on the numerical results obtained using the adaptive multiple-shooting method proposed in this paper. We compare its performance with that of its peers, namely a method based on wavelet-collocation introduced in \cite{EV}, a finite-difference scheme first analyzed
in \cite{AN} and adopted here to solve multi-point SBVPs (see the Appendix for details of its derivation) and a simple-shooting method when it applies.

We have selected three test problems from the literature each exemplifying different characteristics of the solution process. The first problem is a 1-dimensional SBVP with a functional boundary condition and additive noise but the other two are linear 2-dimensional TP-SBVPs, the first with additive and the second with multiplicative noise.
All of the algorithms are implemented in the MATLAB problem-solving environment and executed on a core i5 processor, 2.4GHz, 4GB RAM computer.  \\
\\
\textbf{Test Problem 1}: In this numerical experiment, we try to solve the following 1-dimensional SBVP with functional boundary condition
\begin{equation}\label{T1}
\left\{ {\begin{array}{l}
  \vspace{.2cm}
   {dX(t) = 1\circ dW(t), } \hspace{2cm} 0 \leq t \leq 1, \\
   {\int_{0}^{1}X(s)ds = 0,}  \\
\end{array}} \right.
\end{equation}
and having the exact solution \cite{EV}
 \begin{equation}\label{exact2}
 X(t) =  - \int_{0}^{1}W(s)ds  + W_t.
 \end{equation}

The integral boundary condition in (\ref{T1}) should be discretized (e.g. by the trapezoidal method) into a multi-point boundary condition of the form
$$\frac{\Delta \tau}{2}X(\tau_1)+\sum_{j=2}^{N_s-1}\Delta \tau X(\tau_j)+\frac{\Delta \tau}{2}X(\tau_{N_s})=0,$$
in which $\tau_j=(j-1)\Delta \tau$ for $j=1,\cdots,N_s$ is the $j$-th switching point and $\Delta \tau = \frac{1}{N_s-1}$.
We now place $N_m$ equally-spaced points on the interval $I_i=[\tau_i,\tau_{i+1}]$ which act as the base mesh to integrate the resulting local SDE-IVPs and the global SDE problem. For each realization of the Wiener process (constructed on the base mesh), we solve the system of equations (\ref{NLEQ}) for $\overline{\theta}_{\tau_j}$,~$j=1,2,\cdots,N_s$ and interpolate them by (\ref{interp1}) to arrive at a globally-defined piecewise linear initial approximation to the solution on the whole unit interval.

To find the location of shooting points on $I_i$, we start to synchronously integrate (\ref{SDE-IVP-1}) and (\ref{SDE-IVP-2}) on the base mesh with the schemes described in Section 3.1 to arrive at the first point satisfying (\ref{Stop-Loss}) with $\alpha=0$ and $\beta=2.5$, from where we turn back to $\theta(t)$ and continue the process to reach $\tau_{i+1}$. Similar procedure must be repeated for other intervals to find the set of all optimal shooting points on $[0,1]$.

We are now ready to form and solve the nonlinear system ($\ref{F}$) by a damped-Newton iteration (adopted from \cite{D}) to obtain the optimal starting values at the shooting points (${\bf s}_i$, $i=1,2,\cdots,N_s$ in (\ref{Nonlinear})) and finally solve the original SDE with these initial values by the underlying (full) R3 scheme.

The accuracy of different schemes is measured via the measure $E_{\infty}$ which is an average of the form
$$E_{\infty}:=\frac{1}{M}\sum_{k=1}^{M}E(\omega_k)$$
over a fixed number $M$ of realizations from the maximum grid-wise error
$$E(\omega_k):=\max_{i=1,2,\cdots,N} |\overline{X}(t_i,\omega_k)-X(t_i,\omega_k)|$$
approximating the expected supremum norm
$${\Bbb E}\Big(\|\overline{X}(t)-X(t)\|_{\infty}\Big)$$
in which $\overline{X}(t)$ is our approximation and $X(t)$ is the exact solution. We have used the \texttt{quadl} function in MATLAB to approximate the integral term in the exact solution (\ref{exact2}) which uses the adaptive Lobatto quadrature method.

The results of our computations are depicted in Tables 1 and 2. To compare the accuracy over a single realization ($M=1$), we have provided Table 1 with columns reporting the global error ($E_{\infty}$) for two different methods and a range of grid spacings in the problem domain. The meaning of $N$ in the wavelet-based method is the number of collocation points and in the adaptive multiple-shooting method (or adaptive MSM for short) is the size of base grid used in the integration process. It must be noted that the number of switching points ($N_s$) we have used in the $i$-th row is chosen to be $2^i$ and the number of mid-points ($N_m$) is set accordingly. The superior accuracy of the proposed method (granting one order of magnitude more precision in the results) is obvious from this table. We have observed similar patterns of error behavior over many realizations ($M>>1$) for the adaptive multiple-shooting but due to the unavailability of the data for the other scheme, we have not included them in the Table \ref{Tab1}.
\begin{table}[ht!]\label{Tab1}
\caption{Comparing the accuracy of the wavelet-based and adaptive multiple-shooting methods.}
\centering
\vspace{.5cm}
\begin{tabular}{ccc}
   \hline
   $N$& Wavelet-Collocation & Adaptive MSM \\
   \hline
  $2^{2}$ & 0.2058$\hspace{.5cm}$& 0.0266 \\
  $2^{4}$ & 0.0997$\hspace{.5cm}$ & 0.0036 \\
  $2^{6}$ & 0.0075$\hspace{.5cm}$& 0.0007 \\
  \hline
\end{tabular}
\end{table}

We also have compared the proposed method with a finite-difference scheme in Table 2. Here the errors are reported over $M=500$ realizations (in both methods) and the finite difference equations are set up on all of the $N$ grid-points of the base grid in our adaptive scheme. The column with the heading $N_a$ indicates the average number of shooting points selected by the algorithm. Again we observe a higher accuracy for the adaptive method and a rapid rate of convergence to the exact solution.
\begin{table}[ht!]\label{Tab2}
\caption{Comparing the accuracy of the finite-difference and adaptive multiple-shooting methods.}
\centering
\vspace{.5cm}
\begin{tabular}{cccccc}
  \hline
   $N$&$N_s$&$N_m$&$N_a$& FD Method& Adaptive MSM\\
  \hline
 $2^{5}$&$7$&$4$&11&0.4819&0.0377\\
 $2^{6}$&$10$&$6$&16&0.2728&0.0264\\
 $2^{7}$&$15$&$8$&23&0.2477&0.0164\\
 $2^{8}$&$22$&$12$&36&0.2409&0.0111\\
 $2^{9}$&$32$&$16$&52&0.2287&0.0074\\
  \hline\\
\end{tabular}
\end{table}

\begin{figure}\label{Stab1}
\centering
\includegraphics[angle=0,
width=.5\textwidth]{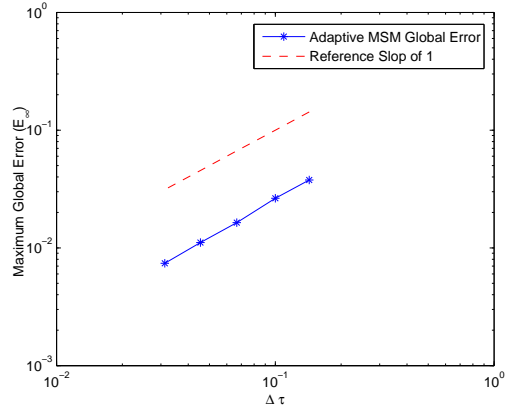}
\caption{Asterisks: strong error measure for the adaptive multiple-shooting method applied to test problem (1). Dashed line: reference slope of 1.}
\end{figure}

In order to investigate the rate of decay of the error (in the strong sense) for the adaptive multiple-shooting method, we have plotted Figure \ref{Stab1} which shows, in a logarithmic scale, the behavior of the global error in terms of increasing the number of switching points. One can observe that the rate of convergence is linear in $\Delta\tau$ and the line of linear regression applied to the data has a slope of $q=1.0126$ with a residual $r=0.0908$. This is a priori anticipated as we have used a method of strong order of convergence one in the integration procedure and a super-linear convergent method in solving the set of nonlinear equations.
\\
\\
\textbf{Test Problem 2}: Here we solve the following 2-dimensional TP-SBVP
\begin{equation*}
\left\{ {\begin{array}{l}
  \vspace{.2cm}
   {dX(t) = (AX(t) + a) dt +(BX(t) + b)\circ dW_t, } \hspace{2cm} 0 \leq t \leq 1, \\
   {H_0X(0) + H_1X(1) = c.}  \\
\end{array}} \right.
\end{equation*}
in which
\begin{equation*}
A = \left[
      \begin{array}{cc}
        0 & 1 \\
        0 & 0 \\
      \end{array}
    \right]
,~a = \left[
      \begin{array}{c}
        0 \\
        c_1 \\
      \end{array}
    \right],~B=\left[
      \begin{array}{cc}
        0 & 0 \\
        0 & 0 \\
      \end{array}
    \right],~
    b = \left[
      \begin{array}{c}
        0 \\
        c_2 \\
      \end{array}
    \right]
    \end{equation*}
    and
 \begin{equation*}
    H_0 = \left[
      \begin{array}{cc}
        1 & 0 \\
        0 & 0 \\
      \end{array}
    \right],~H_1 = \left[
      \begin{array}{cc}
        0 & 0 \\
        1 & 0 \\
      \end{array}
    \right],~
       c = \left[
          \begin{array}{c}
            0 \\
            0 \\
          \end{array}
        \right].
\end{equation*}
Writing $X(t)=[X_1(t),X_2(t)]^T$, it could be shown that $X_1(t)$ solves a second-order SDE having the exact solution
\begin{equation}\label{exact3}
\begin{array}{c}
      \vspace{.2cm}
 X_1(t) = c_1 \frac{t{(t - 1)}}{2} + c_2 (t - 1)\int_0^t {sdW(s)}+t\int_t^1 {(s - 1)dW(s),}
 \end{array}
\end{equation}
and $X_2(t)=\frac{d}{dt}X_1(t)$ (see \cite{AN} for more details).

To obtain the initial trajectory $\theta(t)=[\theta_1(t),\theta_2(t)]^T$ for this test problem and supposing that a realization of $W(1)$ is simulated, we first solve the linear system
$$\Big(H_0+H_1+H_1A+W(1)H_1B\Big)\overline{\theta}_{\tau_1}=c-H_1a-W(1)b$$
for $\overline{\theta}_{\tau_1}$ and then solve another linear system
$$H_1\overline{\theta}_{\tau_2}=c-H_0\overline{\theta}_{\tau_1}$$
for $\overline{\theta}_{\tau_2}$. Now we use linear interpolation to obtain $\theta(t)$ over the whole unit interval. In checking the relation (\ref{Stop-Loss}), we have used the $l_{\infty}$-norm on both sides with $\alpha=2$ and $\beta=1.5$. We also compute the integral terms in the exact solution (\ref{exact3}) by a highly-accurate trapezoidal scheme.

We use the finite-difference and also the simple-shooting methods as two competing approaches to solve this same problem.
Table 3 summarizes our computational results for the case $c_1=1$ and $c_2=1$ averaged over $M=1000$ realizations.
\begin{table}[ht!]\label{Tab3}
\caption{Comparing the accuracy of the finite-difference, simple-shooting and adaptive multiple-shooting methods.}
\centering
\vspace{.5cm}
\begin{tabular}{ccccc}
  \hline
   $N$&$N_a$& FD Method&SS Method& Adaptive MSM\\
  \hline
 $2^{5}$&$9$&$0.0121$&0.0123&1.14e-16\\
 $2^{6}$&$11$&$0.0064$&0.0065&1.70e-16\\
 $2^{7}$&$13$&$0.0032$&0.0032&2.59e-16\\
 $2^{8}$&$15$&$0.0016$&0.0016&4.62e-16\\
 $2^{9}$&$17$&$0.0008$&0.0008&9.62e-16\\
  \hline\\
\end{tabular}
\end{table}
We can observe that while both finite-difference and simple-shooting methods converge uniformly to each other (in terms of accuracy and order of convergence), the adaptive multiple-shooting beats them and gives very accurate results. We also observe a steady growth in the errors as we increase $N$ which could be attributed to the accumulation of round-off errors in the solution process.

In order to show the efficiency of the adaptive method in the weak sense and using the fact that we can compute the expectation of the exact solution and its non-central second moment by the following formulas
\begin{eqnarray}
{\Bbb E}(X_1(t))=c_1\frac{t(t-1)}{2},\quad {\Bbb E}(X_1^2(t))=(3c_1^2+4c_2^2)\frac{t^2(t-1)^2}{12},
\end{eqnarray}
we have approximated these expected values on a range of points in the solution domain by averaging over $M=10000$ realizations of the solution process (computed pointwise) and have compared the results with that of other schemes listed in Table 4. While all methods have a comparable accuracy, the performance of the adaptive method is actually slightly better at all points in the range $[0,1]$.
\begin{table}[ht!]\label{Tab4}
\caption{Comparing the accuracy in the weak-sense of the Heun simple-shooting, finite-difference and adaptive multiple-shooting methods.}
\small
\centering
\vspace{.5cm}
\begin{tabular}{cccccc}
   &Heun Simple-Shooting \cite{AN}& FD Method \cite{AN} &R3 Adaptive MSM & Exact \\
   \hline
   \texttt{t}&  $E(X(t))$ \hspace{.5cm}$E(X^2(t))$ & $E(X(t))$ \hspace{.5cm}$E(X^2(t))$ & $E(X(t))$ \hspace{.5cm}$E(X^2(t))$&$E(X(t))$ \hspace{.5cm}$E(X^2(t))$ \\
   \hline
  0.0  &-0.0000  \hspace{.5cm}0.00000 &-0.0000 \hspace{.5cm} 0.00000 & -0.0000 \hspace{.5cm} 0.0000 &
  -0.0000 \hspace{.5cm} 0.0000 \\
  0.2  & -0.0800  \hspace{.5cm}0.01499&-0.0805 \hspace{.5cm} 0.01497 & -0.0800\hspace{.5cm}  0.0151 & -0.0800\hspace{.5cm}  0.0149\\
  0.4  & -0.1194  \hspace{.5cm}0.03338&-0.1201\hspace{.5cm}  0.03357 & -0.1199\hspace{.5cm}  0.0338 & -0.1200\hspace{.5cm}  0.0336\\
  0.6  & -0.1192  \hspace{.5cm}0.03346&-0.1193\hspace{.5cm}  0.03328 & -0.1202\hspace{.5cm}  0.0338 & -0.1200\hspace{.5cm}  0.0336\\
  0.8  & -0.0793  \hspace{.5cm}0.01486&-0.0791\hspace{.5cm}  0.01472 & -0.0800\hspace{.5cm}  0.0150 & -0.0800\hspace{.5cm}  0.0149 \\
  0.1  & -0.0000  \hspace{.5cm}0.00000 &-0.0000\hspace{.5cm}  0.00000 & -0.0000\hspace{.5cm}  0.0000 & -0.0000\hspace{.5cm}  0.0000 \\
  \hline
\end{tabular}
\end{table}
\\
\textbf{Test Problem 3}: As the last example, we solve the 2-dimensional SDE-BVP system (adopted from \cite{OP}) of the form:
\begin{equation*}
\left\{ {\begin{array}{l}
  \vspace{.2cm}
   {dX(t) = B_1X(t)\circ dW_{1}(t)+B_2X(t)\circ dW_{2}(t), } \quad 0 \leq t \leq 1, \\

   {H_0X(0) + H_1X(1) = c,}  \\
\end{array}} \right.
\end{equation*}
in which
\begin{equation*}
B_1 = \left[
      \begin{array}{cc}
        1 & 1 \\
        0 & 0 \\
      \end{array}
    \right]
,~B_2 = \left[
      \begin{array}{cc}
        0 & 0 \\
        0 & 1 \\
      \end{array}
    \right],
\end{equation*}
and
\begin{equation*}
    H_0 = \left[
      \begin{array}{cc}
        1 & 1 \\
        0 & 0 \\
      \end{array}
    \right],~ H_1 = \left[
      \begin{array}{cc}
        0 & 0 \\
        0 & 1 \\
      \end{array}
    \right],~
    c = \left[
          \begin{array}{c}
            1 \\
            1 \\
          \end{array}
        \right].
\end{equation*}
This equation has an exact solution of the form
\begin{eqnarray}\label{exact4}
X(t)= \left[
        \begin{array}{c}
          e^{W_1(t)}\Big(1-e^{-W_2(1)}+\alpha_{t}^0e^{-W_2(1)}\Big) \\
          e^{W_2(t)}-e^{W_2(1)} \\
        \end{array}
      \right]
\end{eqnarray}
where
\begin{equation}\label{alpha}
\alpha_{t}^s=e^{W_1(t)}\int_s^te^{-W_1(u)}e^{W_2(u)-W_2(s)}dW_1(u).
\end{equation}
Similar to test problem (2), we could obtain the initial trajectory $\theta(t)=[\theta_1(t),\theta_2(t)]^T$ by first simulating a realization from  $W(1)=[W_1(1),W_2(1)]^T$ and then solving the two linear systems
\begin{eqnarray}
\Big(H_0+H_1(I+B_1W_1(1)+B_2W_2(1))\Big)\overline{\theta}_{\tau_1}&=&c,\\
H_1\overline{\theta}_{\tau_2}&=&c-H_0\overline{\theta}_{\tau_1}
\end{eqnarray}
for $\overline{\theta}_{\tau_1}$ and $\overline{\theta}_{\tau_2}$ respectively. Now $\theta(t)$ is computed by linear interpolation and the integration is started to obtain the location of shooting points in the base grid. We use the absolute values of the second components of $\widehat{X}$, $\widetilde{X}$ and $\theta$ in (\ref{Stop-Loss}) with $\alpha=1.5$ and $\beta=2$ and approximate the integral terms in (\ref{exact4}) and (\ref{alpha}) by a sufficiently accurate trapezoidal scheme.

The results of our computations for this test problem are reported in Table 5. For comparison purposes, we have also included the results of applying fixed-step multiple-shooting method in this table. To be fair in the competition, we have selected the number of shooting points in the fixed-step multiple-shooting equal to the average number of adaptive shooting pointes ($N_a$) selected by the adaptive algorithm.

\begin{figure}\label{Stab2}
\centering
\includegraphics[angle=0,
width=.5\textwidth]{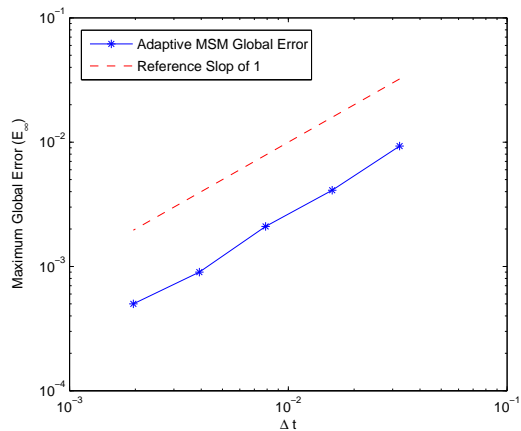}
\caption{Asterisks: strong error measure for the adaptive multiple-shooting method applied to test problem (3). Dashed line: reference slope of 1.}
\end{figure}

In order to examine the strong order of convergence of the adaptive scheme, we have prepared Figure \ref{Stab2} which shows clearly (and in a logarithmic scale) that this order is one. The result of linear regression applied to the data used in the figure gives us a slope of $q=1.0515$ with residual $r=0.1190$ which is acceptable.

\begin{table}[ht!]\label{Tab5}
\caption{Comparing the accuracy of the fixed-step multiple-shooting and adaptive multiple-shooting methods.}
\centering
\vspace{.5cm}
\begin{tabular}{cccccc}
  \hline
   $N$&$N_a$&Fixed MSM Method& Adaptive MSM R3\\
  \hline
 $2^{5}$&$3$&$0.0137$&$0.0093$\\
 $2^{6}$&$4$&$0.0049$&$0.0041$\\
 $2^{7}$&$4$&$0.0025$&$0.0021$\\
 $2^{8}$&$5$&$0.0012$&$0.0009$\\
 $2^{9}$&$5$&$0.0006$&$0.0005$\\
  \hline\\
\end{tabular}
\end{table}
\section{Concluding Remarks}
The numerical solution of boundary value problems in stochastic differential equations is a highly unexplored territory of the SDE world requiring the special attention of the experts in the field to devise methods of high accuracy and efficiency with low computational demand and complexity. We have proposed in this paper, an adaptive multiple-shooting method for general multi-point SBVPs based on a stochastic Runge-Kutta integrator. Although the adaptation criteria is simple and easily implementable in the method, it gives acceptable results in comparison with some other non-adaptive alternatives proposed in the literature. The next step in our research (as explained briefly in Section 3) is to make use of more elaborate stopping criteria in the selection of shooting points and its theortical analysis. We could also incorporate the idea of adptive time-stepping in the integration process itself which we anticipate to improve the accuracy further but needs a theoretical foundation to prove the stability of the overall scheme in a unified manner. Finally, we must also mention the need for introduction of nonlinear test instances into the field which is of great importance for testing and benchmarking purposes in the algorithmic developements expected to be seen in the near future.
\appendix
\section{Finite-Difference Method for Multi-Point SBVPs}
In this appendix, we present a finite-difference scheme for multi-point SBVPs of the form
\begin{equation}\label{SBVP-Ope}
L[X](t) = dW(t)
\end{equation}
in which the operator $L$ is defined by
\begin{equation*}
L = D^n + a_{n-1}(t)D^{n} +\cdots + a_{1}(t)D^{1} + a_{0}(t), \quad D :=\frac{d}{dt},
\end{equation*}
where the coefficients $a_{i}(t)$'s, $i=0,1,\cdots,n-1$ are continuous functions defined on $[0,1]$.
We also append (\ref{SBVP-Ope}) with boundary conditions of the form
\begin{equation}\label{Fun_BC1}
\sum_{j=1}^{N_{s}}\alpha_{ij} X(\tau_{j}) = c_{i}, \hspace{1cm} 1\leq i \leq n,
\end{equation}
defined on some switching points $0\leq\tau_{j}\leq 1,~j=1,2,\cdots,N_s$ (see \cite{AF} for a detailed study of some features of the solution to these problems).

Similar to ordinary differential equations, the SBVP (\ref{SBVP-Ope})-(\ref{Fun_BC1}) can be turned into a first order system
\begin{equation}\label{FOSS}
d{\bf Y(t)} + A(t) {\bf Y(t)} = d{\bf W(t)},
\end{equation}
constrained to satisfy
\begin{equation}\label{Fun_BC2}
\sum_{j=1}^{N_{s}}\alpha_{ij} Y_{n}(\tau_{j}) = c_{i}, \hspace{1cm} 1\leq i \leq n,
\end{equation}
in which ${\bf Y(t)} = (Y_{1}(t),\cdots,Y_{n}(t))$, $Y_{i}(t) = D^{n-i}X(t)$ for $1\leq i \leq n$,
${\bf W}(t) = (W(t),0,\cdots,0)$ and
\begin{equation*}\label{A-Matrix}
A(t) =
\left[
  \begin{array}{ccccc}
    a_{n-1}(t) & a_{n-2}(t) & \cdots & a_{1}(t) & a_{0}(t) \\
    -1 & 0 & \cdots & 0 & 0 \\
    0 & -1 & \cdots & 0 & 0 \\
    \vdots & \vdots & \ddots & \vdots & \vdots \\
    0 & 0 & \cdots & -1 & 0 \\
  \end{array}
\right].
\end{equation*}
To solve (\ref{FOSS})-(\ref{Fun_BC2}) by the finite-difference method, we first construct a base mesh including the switching points of the form
\begin{equation*}\label{Fun_BC}
0 = t_{1}<t_{2}<\cdots<t_{N-1}<t_{N}=1
\end{equation*}
and use an explicit one-step difference scheme on this mesh to arrive at
\begin{equation*}\label{Fun_BC}
{\bf Y}^{j+1}-{\bf Y}^{j} + A^j {\bf Y}^{j} = \Delta{\bf W}^j,
\end{equation*}
in which ${\bf Y}^{j}=[Y_1^j,Y_2^j,\ldots,Y_n^j]^T$ is an approximation to ${\bf Y}(t_j)$, $A^j:=A(t_j)$ and $\Delta{\bf W}^j:={\bf W}(t_{j+1})-{\bf W}(t_j)$.
Simplifying the above relation we obtain
\begin{equation*}\label{Fun_BC}
{\bf Y}^{j+1} + (A^j -I) {\bf Y}^{j} = \Delta{\bf W}^j
\end{equation*}
and arranging them in a sequential manner (into a linear structure) we reach to
\begin{equation}\label{LS}
{\Lambda}{\bf \tilde{Y}} = {\bf w}
\end{equation}
in which
$$ \Lambda_j = A_j - I,\quad
\Lambda = \left[
       \begin{array}{ccccc}
       \Lambda_1 & I &  &   &  \\
       & \Lambda_2 & I &  &  \\
       &  & \ddots &\ddots  &  \\
       &  &  & \Lambda_{N-1} & I \\
            \Phi_1 & \Phi_2  &   \ldots  & \Phi_{N-1} & \Phi_{N} \\
    \end{array}
  \right]_{(N\times n)\times(N\times n)}$$
and
\begin{equation*}
\Phi_j  := \left[ {\begin{array}{*{20}c}
    0& 0 & \ldots & 0 &\alpha_{1j} \\
    0& 0 & \ldots & 0 &\alpha_{2j} \\
   \vdots & &  \ddots& \vdots &\vdots \\
    0& 0 & \ldots & 0 &\alpha_{nj}  \\
\end{array}} \right]=
 \left[
   \begin{array}{c}
     \alpha_{1j} \\
     \alpha_{2j} \\
     \vdots \\
     \alpha_{nj} \\
   \end{array}
 \right]
 \times [0,0,\ldots,0,1]_{1\times n}.
\end{equation*}
Also the vector ${\bf{\tilde Y}}$ has the form
\begin{equation*}
{\bf{\tilde Y}} = [{\bf Y}^1,{\bf Y}^2,\ldots,{\bf Y}^N]^T
\end{equation*}
and
\begin{equation*}
{\bf w} = [\Delta{\bf W}^1,\Delta{\bf W}^2,\ldots,\Delta{\bf W}^N]^T.
\end{equation*}
By solving (\ref{LS}) for each realization of the Wiener process, one obtains the corresponding realization for the solution process on the base mesh which is what we have reported for test problem (1).


\end{document}